\documentclass{article}
\usepackage{amsmath}[1996/11/01]
\usepackage{amssymb,amsthm,amsxtra}

\def\[#1\]{\begin{equation}#1\end{equation}}
\makeatletter
\def\beq{%
   \relax\ifmmode
      \@badmath
   \else
      \ifvmode
         \nointerlineskip
         \makebox[.6\linewidth]%
      \fi
      $$
   \fi
}
\def\eeq{%
   \relax\ifmmode
      \ifinner
         \@badmath
      \else
         $$
      \fi
   \else
      \@badmath
   \fi
   \ignorespaces
}

\def\enddisplaymath{\eeq\global\@ignoretrue}
\makeatother

\newtheorem{thm}{Theorem}
\newtheorem{cor}[thm]{Corollary}
\newtheorem{lem}[thm]{Lemma}
\newtheorem{prop}[thm]{Proposition}

\theoremstyle{definition}

\theoremstyle{remark}
\newtheorem*{rem}{Remark}

\numberwithin{equation}{section}
\numberwithin{thm}{section}

\newcommand{\Q}{\mathbb Q}
\newcommand{\C}{\mathbb C}

\newcommand{\Z}{\mathbb Z}

\renewcommand{\P}{\mathbb P}

\DeclareMathOperator\pf{pf}
\DeclareMathOperator\Hom{Hom}
\DeclareMathOperator\Aut{Aut}

\newcommand{\la}{\langle}
\newcommand{\ra}{\rangle}
\newcommand{\z}{\mathbf{z}}

\newcommand{\II}{\mathord{I\!I}}

\begin{document}

\title{Recurrences for elliptic hypergeometric integrals}
\author{Eric M. Rains\footnote{Department of Mathematics, University of California, Davis}}

\date{April 13, 2005}
\maketitle

\begin{abstract}
In recent work on multivariate elliptic hypergeometric integrals, the
author generalized a conjectural integral formula of van Diejen and
Spiridonov to a ten parameter integral provably invariant under an action
of the Weyl group $E_7$.  In the present note, we consider the action of
the affine Weyl group, or more precisely, the recurrences satisfied by
special cases of the integral.  These are of two flavors: linear
recurrences that hold only up to dimension 6, and three families of
bilinear recurrences that hold in arbitrary dimension, subject to a
condition on the parameters.  As a corollary, we find that a codimension
one special case of the integral is a tau function for the elliptic
Painlev\'e equation.

\end{abstract}


\section{Introduction}

In \cite{xforms}, we studied the following hypergeometric integral
(generalizing the ``Type II'' integral of
\cite{vanDiejenJF/SpiridonovVP:2001}), defined for $|p|,|q|,|t|<1$,
$t_0,\dots,t_7\in \C^*$:
\[
\II^{(n)}_{t;p,q}(t_0,t_1,t_2,t_3,t_4,t_5,t_6,t_7)
:=
\frac{(p;p)^n(q;q)^n}{2^n n!}
\int_{C^n}
\prod_{1\le i<j\le n}
\frac{\Gamma_{p,q}(t z_i^{\pm 1} z_j^{\pm 1})}
     {\Gamma_{p,q}(  z_i^{\pm 1} z_j^{\pm 1})}
\prod_{1\le i\le n}
\frac{\prod_{0\le r\le 7} \Gamma_{p,q}(t_r z_i^{\pm 1})}
     {\Gamma_{p,q}(z_i^{\pm 2})}
\frac{dz_i}{2\pi\sqrt{-1}z_i},
\]
where $(p;p)=\prod_{1\le i} (1-p^i)$, $\Gamma_{p,q}$ is the elliptic Gamma
function (see below), and $C$ is a suitable choice of contour (which may be
taken to be the unit circle when all parameters are inside the unit circle).
We found that if the parameters satisfied the following ``balancing''
condition:
\[
t^{2n-2}t_0t_1t_2t_3t_4t_5t_6t_7=p^2q^2,
\]
then the integral satisfied a certain transformation which, together with
the obvious permutation symmetry of the arguments, generated an action of
the Weyl group $W(E_7)$.  More precisely, assuming balanced parameters, the
renormalized integral
\[
\tilde{\II}^{(n)}_{t;p,q}(t_0,t_1,\dots,t_7)
:=
\!\!\prod_{0\le r<s\le 7} \!\!\Gamma^+_{t,p,q}(tt_rt_s)\ 
\II^{(n)}_{t;p,q}(t^{1/2}t_0,t^{1/2}t_1,\dots,t^{1/2}t_7)
\]
is invariant under this action, which we now explain.

We first observe that we can view the above integral (given the balancing
condition) as a function on an algebraic torus, the maximal torus
$\Hom(\Lambda_{E_8},\C^*)$ of the complex Lie group $E_8$ (where
$\Lambda_{E_8}$ is the root lattice).  Indeed, we first observe that the
integral is invariant under the symmetry
\[
(t_0,t_1,\dots,t_7)\mapsto (-t_0,-t_1,\dots,-t_7),
\]
simply by negating the $z$ variables; as a result, it is only a function of
the pairwise products and ratios $t_r^{\pm 1} t_s^{\pm 1}$.  In other
words, it is a function on the maximal torus $\Hom(\Lambda_{D_8},\C^*)$.
But the balancing condition forces a choice of square root
\[
\sqrt{t_0t_1\dots t_7} = pq/t^{n+1}
\]
and thus the parameters in fact determine a homomorphism from the lattice
$\Lambda_{E_8}$ to $\C^*$, mapping $\omega:=(1/2,1/2,\dots,1/2)$ to
$pq/t^{n+1}$.

If $\phi:\Lambda_{E_8}\to C^*$ is a homomorphism such that
\[
\frac{pq}{t\phi(\omega)}=t^n
\]
for some (uniquely determined) integer $n$, we define
\[
\tilde\II_{t;p,q}(\phi)
\]
as follows.  If $n<0$, then
\[
\tilde\II_{t;p,q}(\phi)=0;
\]
otherwise, we set
\[
\tilde\II_{t;p,q}(\phi)
=
\tilde\II^{(n)}_{t;p,q}(\phi(e_0),\dots,\phi(e_7))
\]
where $e_0,\dots,e_7$ are the coordinate vectors and we have chosen an
extension of $\phi$ to $\Lambda_{D_8}^*$ (which as remarked above does not
affect the value of the integral).

\begin{thm}\cite{xforms}
Suppose $\phi\in\Hom(\Lambda_{E_8},\C^*)$ Then for any element $g\in
W(E_8)=\Aut(\Lambda_{E_8})$ such that
\[
\la \omega,g\omega\ra\in \{1,2\},
\]
we have
\[
\tilde\II_{t;p,q}(\phi)
=
\tilde\II_{t;p,q}(g^*\phi)
\]
whenever
\[
\frac{pq}{t\phi(\omega)},
\frac{pq}{t\phi(g\omega)}
\in t^\Z,
\]
so that both sides are defined.
\end{thm}

Note that if $\la \omega,g\omega\ra=2$, then $g\omega=\omega$.  In other
words, $g$ is in the stabilizer $W(E_7)$ of $\omega$, and the statement
becomes that $\tilde\II_{t;p,q}(\phi)$ is invariant under $W(E_7)$ whenever
it is defined.

In addition to the natural action of the finite Weyl group $W(E_8)$ on
$\Hom(\Lambda_{E_8},\C^*)$, there is a nearly natural action of the affine
Weyl group.  To be precise, if $v\in \Lambda_{E_8}$, we define a
shift operator $\tau_v$ by
\[
(\tau_v(\phi))(w) = \phi(w) q^{\la v,w\ra},
\]
for all $\phi\in \Hom(\Lambda_{E_8},\C^*)$, $w\in \Lambda_{E_8}$.  (It will
be notationally convenient to extend this definition to
$v\in \Lambda_{E_8}\otimes \Q$ by fixing a consistent family of $m$th roots
of $q$.)  The price of enlarging the group is that we no longer have
invariance; instead, the most we can expect is that $\tilde\II_{t;p,q}$
should satisfy recurrences with respect to different shifts.

The purpose of the present note is to show that in certain special cases,
such recurrences do indeed arise.  These come in two main flavors.  The
first set of recurrences arises from the observation that certain shifts
(by coordinate vectors, say) have the effect of multiplying the integrand
by a relatively simple function; in low dimensions ($n\le 6$), these
functions must be linearly dependent, and thus give rise to a linear
recurrence.

The other set of recurrences are somewhat more subtle.  The above integral
can be viewed as a generalization of the Selberg integral, which suggests
that the speical cases $t\in \{q^{1/2},q,q^2\}$ should be particularly
nice.  Indeed, it turns out that in those cases the integral can be
expressed (in many ways) as a determinant or pfaffian of one- or
two-dimensional integrals.  In particular, we can arrange for several
minors of said determinant/pfaffian to themselves be special cases of our
integral, with the result that the Pl\"ucker relations give rise to
recurrences of our integral.  Since the Pl\"ucker relations are bilinear,
the resulting recurrences are also bilinear; for $t=q$ (the determinantal
case), we obtain a three-term bilinear recurrence, while for
$t=q^{1/2},q^2$ (pfaffian cases), we obtain a four-term bilinear
recurrence.  The significance of these recurrences is perhaps underscored
by the fact that the recurrence for $t=q$ has arisen in the theory of
Sakai's elliptic Painlev\'e equation
\cite{SakaiH:2001,KajiwaraK/MasudaT/NoumiM/OhtaY/YamadaY:2004}.

The plan of the paper is as follows.  After defining some notation for
generalized $q$-symbols and theta functions, we proceed in section 2 to
prove some theta function identities needed in the derivation of our
recurrences.  In section 3, we use these to give the aforementioned linear
recurrences in low dimensions.  Section 4 describes a general setting in
which Pl\"ucker relations give rise to bilinear relations of integrals,
which is then specialized in section 5 to give our bilinear Painlev\'e-type
recurrences.

\noindent{\bf Acknowledgements}

The author would like to thank M. Noumi and K. Takasaki for organizing
such an interesting workshop, as well as M. Adler and P. van Moerbeke
for helpful conversations on Fay identities and tau functions.  The
author's work was partially supported by NSF Grant No. DMS-0401387.

\noindent{\bf Notation}

Aside from the integral itself, most of the functions that appear in the
sequel are most simply expressed as infinite products; as a result, we will
need a shorthand notation for certain such products.  Here $p$, $q$, $t$
are complex numbers inside the open unit disc.
\begin{align}
\theta_p(x) &:= \prod_{0\le k} (1-p^{k+1}/x)(1-p^k x)\\
\Gamma_{p,q}(x) &:= \prod_{0\le j,k} (1-p^{j+1}q^{k+1}/x)(1-p^j q^k x)^{-1}\\
\Gamma^+_{p,q,t}(x) &:= \prod_{0\le i,j,k} (1-p^{i+1}q^{j+1}t^{k+1}/x)(1-p^i q^j t^k x).
\end{align}
The first function is simply a version of Jacobi's theta function, while
the second function is Ruijsenaars' elliptic Gamma
function \cite{RuijsenaarsSNM:1997}.  As these are generalized $q$-symbols
(indeed, $\Gamma_{0,q}(x)^{-1}$ is precisely the usual $q$-symbol), we take
the standard convention that the presence of multiple arguments indicates a
product; thus, for instance, in the above integral,
\[
\Gamma_{p,q}(z_i^{\pm 1}z_j^{\pm 1}) =
\Gamma_{p,q}(z_iz_j)
\Gamma_{p,q}(z_i/z_j)
\Gamma_{p,q}(z_j/z_i)
\Gamma_{p,q}(1/z_iz_j).
\]
The main properties of these functions are reflection symmetry:
\begin{align}
\theta_p(p/x)&=\theta_p(x)\\
\Gamma_{p,q}(pq/x)&=\Gamma_{p,q}(x)^{-1}\\
\Gamma^+_{p,q,t}(pqt/x)&=\Gamma_{p,q,t}(x),
\end{align}
and a functional equation:
\begin{align}
\theta_p(px) &= \frac{1-1/x}{1-x} \theta_p(x) = -x^{-1}\theta_p(x)\\
\Gamma_{p,q}(qx) &= \theta_p(x)\Gamma_{p,q}(x)\\
\Gamma^+_{p,q,t}(tx) &= \Gamma_{p,q}(x)\Gamma^+_{p,q,t}(x),
\end{align}
with similar identities following by the symmetry of $\Gamma_{p,q}$ and
$\Gamma^+_{p,q,t}$ in the parameters.

We recall that a ($p$-)theta function (in multiplicative notation) is a
holomorphic function $f(x)$ on $\C^*$ such that
\[
f(px) = C (-x)^{-m} f(x)
\]
for some constant $C$ (the multiplier), and some integer $m$ (the
degree).  The canonical example of this is the function $\theta_p(x/a)$;
indeed, any $p$-theta function is proportional to a function of the form
\[
x^k \prod_{1\le i\le m} \theta_p(x/a_i),
\]
with multiplier
\[
p^k \prod_{1\le i\le m} a_i,
\]
and thus the multiplier of a theta function is determined up to powers of
$p$ by its zeros.  A meromorphic theta function is a ratio of holomorphic
theta functions.

Similarly, a $BC_n$-symmetric theta function of degree $m$ is defined to be
a function on $(\C^*)^n$ invariant under permutations and inversions of its
variables, and such that as a function of each variable it is a theta
function of degree $2m$ with multiplier $p^{-m}$.  Since the quotient of
the elliptic curve $\C^*/\la p\ra$ by $x\mapsto 1/x$ is a projective line,
it follows that the space of $BC_1$-symmetric theta functions of degree $m$
is $m+1$-dimensional.

\section{Theta function relations}

Define a function $\psi_p(x,y)$ on $\C^*\times \C^*$ as follows:
\[
\psi_p(x,y) = x^{-1}\theta_p(xy)\theta_p(x/y).
\]
This is readily seen to satisfy the relations
\[
\psi_p(x,y)=\psi_p(x,1/y)=-\psi_p(y,x)
\]
and
\[
\psi_p(x,py)=(py^2)^{-1}\psi_p(x,y).
\]
Somewhat less trivial is the following:

\begin{lem}
For $x,y,z,w\in \C^*$,
\[
\psi_p(x,y)\psi_p(z,w)-\psi_p(x,z)\psi_p(y,w)+\psi_p(x,w)\psi_p(y,z)
=
0.
\]
\end{lem}

\begin{proof}
Consider the skew-symmetric $4\times 4$ matrix
\[
A=
\begin{pmatrix}
\psi_p(x,x)&\psi_p(x,y)&\psi_p(x,z)&\psi_p(x,w)\\
\psi_p(y,x)&\psi_p(y,y)&\psi_p(y,z)&\psi_p(y,w)\\
\psi_p(z,x)&\psi_p(z,y)&\psi_p(z,z)&\psi_p(z,w)\\
\psi_p(w,x)&\psi_p(w,y)&\psi_p(w,z)&\psi_p(w,w)
\end{pmatrix}
\]
The functions $\psi_p(x,\_)$, $\psi_p(y,\_)$, $\psi_p(z,\_)$,
$\psi_p(w,\_)$ all lie in the 2-dimensional space of $BC_1$-symmetric theta
functions of degree 1, and thus any three of them satisfy a linear
relation.  In particular, it follows that the matrix $A$ has rank at most
2, and thus has pfaffian 0; this is precisely the desired identity.
\end{proof}

\begin{rem} This, of course, is simply the addition law for elliptic theta
  functions in disguise.
\end{rem}

\begin{prop}\label{prop:ell_CauchyII}
We have the following Cauchy-type determinant:
\[
\det_{1\le i,j\le n}(\frac{1}{\psi_p(x_i,y_j)})
=
(-1)^{n(n-1)/2}
\frac{\prod_{1\le i<j\le n} \psi_p(x_i,x_j)\psi_p(y_i,y_j)}
     {\prod_{1\le i,j\le n} \psi_p(x_i,y_j)}
\]
\end{prop}

\begin{proof}
From the lemma, we can write
\begin{align}
\frac{1}{\psi_p(x_i,y_j)}
&=
\frac{\psi_p(z,w)}
     {\psi_p(z,x_i)\psi_p(w,y_j)-\psi_p(z,y_j)\psi_p(w,x_i)}\\
&=
\frac{\psi_p(z,w)}{\psi_p(w,x_i)\psi_p(w,y_j)}
\frac{1}
     {(\psi_p(z,x_i)/\psi_p(w,x_i))-(\psi_p(z,y_j)/\psi_p(w,y_j))}
\end{align}
for arbitrary $z$, $w$.  The result thus follows immediately from the
usual Cauchy determinant.
\end{proof}

\begin{rem}
That this identity is a special case of the usual Cauchy determinant is no
accident: any function $\psi$ satisfying the above identity can be written
in the form
\[
\psi(x,y)
=
\frac{\psi(z,x)\psi(w,y)-\psi(z,y)\psi(w,x)}
     {\psi(z,w)}
\]
using the $n=2$ instance of the identity.
\end{rem}

\begin{cor}\label{cor:BCn_parfrac}
For generic $x_1$,\dots,$x_{n+2}$, $y_1$,\dots,$y_n$,
\[
\sum_{1\le k\le n+2}
\frac{\prod_{1\le j\le n} \psi_p(x_k,y_j)}
     {\prod_{i\ne k} \psi_p(x_k,x_i)}
=
0.
\]
\end{cor}

\begin{proof}
Expand the $n+1$ dimensional instance of the above determinant along the
last row, set $y_{n+1}=x_{n+2}$, then simplify.  Alternatively, observe
that some such relation must hold by dimensionality, and deduce the
constants by setting $y_j=x_k$ for various choices of $j$, $k$.
\end{proof}

Fix $p\in \C$, $t\in \C^*$, and define for $u_0$, $u_1$, $u_2$, $u_3$,
$u_4\in \C^*$ a function $g^{(n)}_{u_0,u_1,u_2,u_3,u_4}$ on $(\C^*)^n$ by
\[
g^{(n)}_{u_0,u_1,u_2,u_3,u_4}(\dots z_i\dots)
=
\prod_{1\le i\le n} (1+R(z_i))
\frac{\prod_{0\le r\le 4} \theta_p(u_r
  z_i)\theta_p(z_i/t^{n-1}u_0u_1u_2u_3u_4)}
     {z_i^2\theta_p(z_i^2)}
\prod_{1\le i<j\le n} \frac{\theta_p(tz_iz_j)}{\theta_p(z_iz_j)},
\]
where $R(z_i)$ is the operator $z_i\mapsto 1/z_i$.  We also define a
function
\[
f^{(n)}_{u_0}(\dots z_i\dots)
=
\prod_{1\le i\le n} \theta_p(u_0 z_i,u_0/z_i).
\]
The following lemma shows that this is a special case of the first family.

\begin{lem}\label{lem:fasg}
We have the identity
\[
g^{(n)}_{u_0,u_1,u_2,u_3,1/u_0}(\dots z_i\dots)
=
f^{(n)}_{u_0}(\dots z_i\dots)
\prod_{1\le i\le n}
\frac{\theta_p(t^{n-i}u_1u_2,t^{n-i}u_1u_3,t^{n-i}u_2u_3)}
     {t^{n-1}u_0u_1u_2u_3}.
\]
\end{lem}

\begin{proof}
If we divide both sides by $f^{(n)}_{u_0}(\dots z_i\dots)$, the result is
simply Lemma 6.2 of \cite{xforms}.
\end{proof}

\begin{thm}\label{thm:grels}
For $u_0,u_1,u_2,u_3,v_0,\dots,v_{n+1}\in \C^*$, and $\z\in (\C^*)^n$,
\[
\sum_{0\le i\le n+1}
\frac{g^{(n)}_{u_0,u_1,u_2,u_3,v_i}(\z)}
     {\prod_{r\ne i} v_r^{-1} \theta_p(v_r/v_i,t^{n-1}u_0u_1u_2u_3v_iv_r)}
=
0.
\]
\end{thm}

\begin{proof}
If we pull the sum inside the symmetrization operation, we find that the
result would follow from the identity
\[
\sum_{0\le k\le n+1}
\frac{\prod_{1\le i\le n} z_i^{-1}\theta_p(v_kz_i,z_i/t^{n-1}u_0u_1u_2u_3v_k)}
     {\prod_{r\ne k} v_r^{-1} \theta_p(v_r/v_k,t^{n-1}u_0u_1u_2u_3v_kv_r)}
=
0.
\]
But this is the special case of Corollary \ref{cor:BCn_parfrac} with
\[
x_k = v_{k-1}\sqrt{t^{n-1}u_0u_1u_2u_3}\quad
y_k = \frac{z_k}{\sqrt{t^{n-1}u_0u_1u_2u_3}}.
\]
\end{proof}

If we set one of the variables in $g^{(n)}_{u_0u_1u_2u_3u_4}$ equal to
$u_0$, half of the terms vanish, and we thus find
\[
g^{(n)}_{u_0u_1u_2u_3u_4}(u_0,\z)
=
\frac{\theta_p(1/t^{n-1}u_1u_2u_3u_4)\prod_{1\le r\le 4} \theta_p(u_0 u_r)}
     {u_0^2}
g^{(n-1)}_{t u_0,u_1,u_2,u_3,u_4}(\z).
\]
This in some cases allows us to deduce relations between these functions.
We concentrate on the case $n=4$, as this seems to be the primary source of
identities between functions $f^{(n)}$ and $g^{(n)}$ not contained in
Corollary \ref{cor:BCn_parfrac} or Theorem \ref{thm:grels}; all other
such identities we have been able to find are obtained by specializing the
variables.

\begin{prop}\label{prop:gtof}
For any parameters $u_0$, $u_1$, $u_2$, $u_3$, $u_4\in\C$, $\z\in \C^{4}$,
\[
g^{(4)}_{u_0,u_1,u_2,u_3,u_4}(\z)
=
\prod_{0\le i<j\le 4} \theta_p(u_iu_j,t u_iu_j)
\sum_{0\le r\le 4}
\prod_{i\ne r}
\frac{\theta_p(u_i/t^3u_0u_1u_2u_3u_4)}
     {u_i^2\theta_p(u_r/u_i,u_ru_i,t u_ru_i)}
f^{(4)}_{u_r}(\z)
\]
\end{prop}

\begin{proof}
The five functions $f^{(4)}_{u_r}(\z)$ span the space of $BC_4$-symmetric
theta functions of degree 1, so it remains only determine the coefficients
of the expansion.  If we evaluate $g^{(4)}_{u_0,u_1,u_2,u_3,u_4}(\z)$ at the
point $\z=(u_1,u_2,u_3,u_4)$, only the $f^{(4)}_{u_0}(\z)$ term survives,
and we thus can solve for its coefficient; the other coefficients are
symmetrical.
\end{proof}

There is a sort of inverse to the above expansion, expressing
$f^{(4)}_{u_0}(\z)$ in terms of the five functions
\[
g^{(4)}_{u_0,u_2,u_3,u_4,u_5}(\z),
g^{(4)}_{u_0,u_1,u_3,u_4,u_5}(\z),\dots,
g^{(4)}_{u_0,u_1,u_2,u_3,u_4}(\z).\notag
\]

\begin{prop}\label{prop:ftog}
For any parameters $u_0$, $u_1$, $u_2$, $u_3$, $u_4$, $u_5\in\C$, $\z\in
\C^{4}$,
\[
f^{(4)}_{u_0}(\z)
=
\frac{
\prod_{1\le i\le 5} \theta_p(u_0 u_i/t^3 U)
}{\prod_{1\le i<j\le 5} \theta_p(u_iu_j,t u_iu_j)}
\sum_{1\le r\le 5}
\frac{g_{u_0,\dots,\widehat{u_r},\dots}(\z)}
     {\theta(u_0u_i/t^3 U)}
\prod_{0\le i\le 3} \frac{\theta_p(t^i u_0u_r)}{\theta_p(t^{i-6}/U)}
\prod_{1\le i\ne r} \frac{u_i^2\theta_p(u_iu_r,t u_iu_r)}{\theta_p(u_r/u_i)},
\]
where $U = u_0u_1u_2u_3u_4u_5$.
\end{prop}

\begin{proof}
From Proposition \ref{prop:gtof} above, we obtain six different identities
expressing the six functions $g^{(4)}_{\dots,\widehat{u_r},\dots}(\z)$ in
terms of the six functions $f^{(4)}_{u_r}(\z)$.  It turns out, in fact,
that up to rescaling of rows and columns, the resulting $6\times 6$ matrix
is antisymmetric, and thus the inverse matrix can be expressed via
pfaffians.  The closed forms for the desired pfaffians can be obtained via
the special case $a=1$, $b=t$, $c=1/t^3U$ of the following identity.
\end{proof}

\begin{thm}\cite{OkadaS:2004}
For arbitrary parameters $u_0,\dots,u_{2n-1},a,b,c\in C$, we have
\[
\pf_{0\le i,j<2n}
\left(\frac{u_j\theta_p(u_i/u_j,a u_iu_j,b u_iu_j)}{\theta_p(c u_iu_j)}\right)
=
c^{n(n-1)}
\theta_p(a/c,b/c)^{n-1}
\theta_p(ac^{n-1}U,bc^{n-1}U)
\prod_{0\le i<j<2n} \frac{u_j\theta_p(u_i/u_j)}{\theta_p(c u_iu_j)},
\]
where $U=\prod_{0\le i<2n} u_i$
\end{thm}

\begin{proof}
We first consider both sides as functions in $a$; we find that they are
both theta functions with the same multiplier.  Moreover, if we set $a=c$
on the left, we obtain a matrix of rank 2, and thus the pfaffian must have
a zero of order $n-1$ at that point.  This accounts for all but one zero in
a fundamental region, and the remaining zero can be determined from the
multiplier.  Arguing similarly for $b$, we conclude that the left-hand side
is a multiple of
\[
\theta_p(a/c,b/c)^{n-1}
\theta_p(ac^{n-1}U,bc^{n-1}U).
\]
Since the pfaffian also vanishes whenever $u_i=u_j$, and has at most simple
poles at points with $cu_iu_j=1$, it follows that the ratio of the two
sides is in fact constant.  The value of this constant can then be
determined from the asymptotics as $u_{2i}\to 1/cu_{2i-1}$.
\end{proof}

Associated to this is the following analogue of Corollary
\ref{cor:BCn_parfrac}.

\begin{cor}
For arbitrary parameters $a,b\in \C$, $u\in \C^{n+1}$,
\[
\sum_{0\le r\le n}
\theta_p(a u_r,b u_r,a U/u_r,b U/u_r)
\prod_{0\le i\le n;i\ne r} \frac{\theta_p(u_iu_r)}{u_i\theta_p(u_r/u_i)}
=
\delta_{\text{$n$ even}}
\theta_p(a,b,aU,bU),
\]
where $U = \prod_{0\le r\le n} u_r$.
\end{cor}

\begin{proof}
The case $n$ even can be obtained by setting $c=u_{2n-1}=1$ and expanding the
pfaffian along the last row; the case $n$ odd then follows by setting
$u_{2n-2}=1$.
\end{proof}

Similarly, the fact that the pfaffians are nice gives rise to a relation
between the functions
\[
f^{(4)}_{u_0}(\z),
f^{(4)}_{u_1}(\z),
f^{(4)}_{u_2}(\z),
g^{(4)}_{u_0u_1u_2u_3u_4}(\z),
g^{(4)}_{u_0u_1u_2u_3u_5}(\z),
g^{(4)}_{u_0u_1u_2u_4u_5}(\z).
\]
If we also use the relation between $f^{(4)}_{u_r}(\z)$, $0\le r\le 5$
coming from Corollary \ref{cor:BCn_parfrac}, then we can obtain similar
relations involving $4$, $2$, or $0$ of the $f$ functions; we omit the
details.

Finally, we will also need the following pfaffian identity.

\begin{thm}\label{thm:ell_CauchypfaffII}
We have the pfaffian
\[
\pf_{1\le i<j\le 2n}(
\frac{z_i^{-1}\theta(z_iz_j^{\pm 1};t^2)}
     {\theta(t z_iz_j^{\pm 1};t^2)}
)
=
\frac{t^{n(n-1)}
\prod_{1\le i<j\le 2n}
z_i^{-1}\theta(z_iz_j^{\pm 1};t^2)}
{\theta(t z_iz_j^{\pm 1};t^2)}
\]
\end{thm}

\begin{proof}
Both sides are $BC_n$-antisymmetric abelian functions with the same polar
divisor, and are thus proportional.  Multiplying both sides by
\[
\prod_{1\le i\le n}
\theta(t z_{2i-1}/z_{2i};t^2)
\]
and taking the limit $z_{2i}\to t z_{2i-1}$ shows that the constant is 1.
\end{proof}

\section{Recurrences in low dimensions}

We can obtain recurrences for low-dimensional instances of our integral by
observing that there are two ways in which shifting the parameters
corresponds to multiplying the integrand by a degree 1 theta function.
If we multiply $t_r$ by q, this simply multiplies the integrand by
\[
f^{(n)}_{t_r}(\dots z_i\dots)
=
\prod_{1\le i\le n} \theta_p(t_r z_i^{\pm 1} )
=
t_r^n \prod_{1\le i\le n} \psi_p(t_r,z_i).
\]
Somewhat more subtly, if $t^{2n-2}t_0t_1t_2t_3t_4t_5t_6t_7=p^2q$,
multiplying the integrand by
\[
\prod_{0\le i<n} 
\frac{q^2p^3}
     {(t_5t_6t_7)^2 t^{2n-2}\theta_p(pq/t^it_5t_6,pq/t^it_5t_7,pq/t^it_6t_7)}
g^{(n)}_{t_0,t_1,t_2,t_3,t_4}(\dots z_i\dots)
\]
simply has the effect of multiplying $t_0$ through $t_4$ by $\sqrt{q}$ and
dividing $t_5$, $t_6$, $t_7$ by $\sqrt{q}$.  Indeed, this follows
immediately by an adjointness argument as in the second proof of Theorem
6.1 of \cite{xforms}.

As a result, any linear dependence between the 8 functions $f^{(n)}$ and
the 56 functions $g^{(n)}$ gives rise to a relation of integrals.  Thus in
principle we would obtain recurrences all the way up to dimension 62 (since
the space of degree 1 $BC_n$-symmetric theta functions has dimension
$n+1$); in practice, however, the coefficients of such relations do not
appear to have nice closed forms in general.  There is, however, one
special case in which the coefficients are nice.  Each function corresponds
to a vector (by which it shifts the parameters); if the difference of any
two such vectors in the collection is a root of $E_7$, the corresponding
relation has nice coefficients.  This, however, greatly reduces the
possible number of theta functions in the relation, with the result that we
only obtain recurrences for $n\le 6$.

The simplest case is the linear relations between the functions $f^{(n)}$
from Corollary \ref{cor:BCn_parfrac}, which gives the following recurrence.

\begin{thm}\label{thm:recurs_ld}
For $1\le n\le 6$, let $t_0$,\dots, $t_7$, $t$, $p$, $q$ be parameters such
that $|p|,|q|,|t|<1$.  Then
\[
\sum_{0\le i\le n+1}
\frac{t_i \II^{(n)}_{t;p,q}(t_0,\dots,q t_i,\dots,t_7)}
     {\prod_{0\le j\le n+1;j\ne i} \theta_p(t_i t_j^{\pm 1})}
=
0.
\]
\end{thm}

Another source of such recurrences is Theorem \ref{thm:grels}, especially in
combination with Lemma \ref{lem:fasg}.  The upshot is that we obtain
(relatively) nice relations between any $n+2$ of the 8 functions
\[
f^{(n)}_{t_0},
f^{(n)}_{t_1},
f^{(n)}_{t_2},
f^{(n)}_{t_3},
g^{(n)}_{t_0,t_1,t_2,t_3,t_4},
g^{(n)}_{t_0,t_1,t_2,t_3,t_5},
g^{(n)}_{t_0,t_1,t_2,t_3,t_6},
g^{(n)}_{t_0,t_1,t_2,t_3,t_7};
\]
we simply apply Theorem \ref{thm:grels} with $u_i=t_i$;
$v_0,\dots,v_{n+1}\in \{1/t_0,1/t_1,1/t_2,1/t_3,t_4,t_5,t_6,t_7\}$.  The
coefficients of the resulting relations are, unfortunately, rather
complicated (albeit products of theta functions).  In fact, the resulting
recurrences are simply images of the recurrence of Theorem
\ref{thm:recurs_ld} under the action of the Weyl group $W(E_7)$ (assuming,
of course, that $t^{2n-2}t_0t_1t_2t_3t_4t_5t_6t_7=p^2q$, so that $W(E_7)$
actually does act).  Ideally, we would prefer to give a manifestly
$W(E_7)$-invariant description of the recurrences; in the absence of such a
description, we leave the details to the reader, rather than list all of
the superficially different recurrences arising in this way.

Another $W(E_7)$-orbit of recurrences arises from Propositions \ref{prop:gtof}
and \ref{prop:ftog} above.  For instance, from Proposition \ref{prop:gtof},
we obtain the following.

\begin{thm}
For $n=4$, let $t_0$,\dots, $t_7$, $t$, $p$, $q$ be parameters such
that $|p|,|q|,|t|<1$, $t^6 t_0t_1t_2t_3t_4t_5t_6t_7 = p^2q$.  Then
\begin{align}
&\II^{(4)}_{t;p,q}(q^{1/2}t_0,q^{1/2}t_1,q^{1/2}t_2,q^{1/2}t_3,q^{1/2} t_4,q^{-1/2}t_5,q^{-1/2}t_6,q^{-1/2}t_7)\\
&=
\frac{t^{24} (t_0t_1\dots t_4)^8 \prod_{0\le i<j\le 4} \theta_p(t_it_j,t t_it_j)}
     {p^4 \prod_{0\le i<4} \theta_p(pq/t^it_5t_6,pq/t^it_5t_7,pq/t^it_6t_7)}
\sum_{0\le r\le 4}
\prod_{0\le i\le 4;i\ne r}
\frac{\theta_p(t_i/t^3t_0t_1t_2t_3t_4)}
     {t_i^2\theta_p(t_r/t_i,t_rt_i,t t_rt_i)}
\II^{(4)}_{t;p,q}(t_0,\dots,qt_r,\dots,t_7)
\notag
\end{align}
\end{thm}

In fact, the $W(E_7)$-images of this identity and those of Theorem
\ref{thm:recurs_ld} include every linear recurrence in which the
differences of any two shifts is a root of $E_7$.

\section{Generalized Fay identities}

Suppose $\psi(x,y)$, $\psi'(x,y)$ are antisymmetric measurable functions on
$X^2$ for some space $X$ that satisfy the identity of Proposition
\ref{prop:ell_CauchyII}; for instance, $\psi(x,y)=\psi_p(x,y)$.  There is a
natural family of multidimensional integrals attached to these functions in
such a way that the Pl\"ucker relations between minors of a matrix
translate into bilinear identities satisfied by integrals.

We define, for any measure $\mu$
\[
\tau^{(n)}(\mu;\psi,\psi')
=
\frac{1}{n!}
\int_{X^n}
\prod_{1\le i<j\le n} \psi(x_i,x_j)
\prod_{1\le i<j\le n} \psi'(x_i,x_j)
\prod_{1\le i\le n} \mu(dx_i),
\]
assuming this integral converges.  In addition, for notational convenience,
we define
\[
\tau^{(n)}(\mu[a_1,\dots,a_k][b_1,\dots,b_l]';\psi,\psi')
:=
\prod_{1\le i<j\le k} \psi(a_i,a_j)
\prod_{1\le i<j\le l} \psi'(b_i,b_j)
\tau^{(n)}(\prod_{1\le i\le k}[a_i]\prod_{1\le i\le
  l}[b_i]'\mu;\psi,\psi')
\]
where $[a_i](x)=1/\psi(a_i,x)$, $[b_i]'(x)=1/\psi'(b_i,x)$.
Since we will for the most part be fixing $\psi,\psi'$, we will suppress
them from the notation when no confusion will result.

Using the fact that $\psi$, $\psi'$ satisfy Cauchy-type identities, we find
that $\tau^{(n)}$ is the integral of a product of two determinants, and
thus by the integral analogue of Cauchy-Binet, is itself a determinant of
univariate integrals, and can be written as such a determinant in many
different ways.

\begin{thm}
Assuming all integrals are defined,
\[
\tau^{(n)}(\mu[a_1,\dots,a_k][b_1,\dots,b_l]')
=
\det_{1\le i,j\le n} \tau^{(1)}(\mu[a_i][b_j]')
.
\]
\end{thm}

\begin{proof}
We have
\[
\prod_{1\le i<j\le n} \psi(x_i,x_j)
=
(-1)^{n(n-1)/2}
\frac{\prod_{1\le i,j\le n} \psi(a_i,x_j)}
     {\prod_{1\le i<j\le n} \psi(a_i,a_j)}
\det_{1\le i,j\le n}(\frac{1}{\psi(a_i,x_j)})
\]
and similarly
\[
\prod_{1\le i<j\le n} \psi'(x_i,x_j)
=
(-1)^{n(n-1)/2}
\frac{\prod_{1\le i,j\le n} \psi'(b_i,x_j)}
     {\prod_{1\le i<j\le n} \psi'(b_i,b_j)}
\det_{1\le i,j\le n}(\frac{1}{\psi'(b_i,x_j)})
\]
and thus
\begin{align}
\tau^{(n)}(\mu[a_1,\dots,a_n][b_1,\dots,b_n]')
&=
\int_{X^n}
\det_{1\le i,j\le n}(\frac{1}{\psi(a_i,x_j)})
\det_{1\le i,j\le n}(\frac{1}{\psi'(b_i,x_j)})
\mu(dx_i)\\
&=
\det_{1\le i,j\le n}
\int_{X}
\frac{1}
{\psi(a_i,x)\psi'(b_j,x)}
\mu(dx).
\end{align}
\end{proof}

Now, any minor of the above determinant is itself a determinant of the same
form.  As a consequence, any polynomial equation satisfied by minors of a
general matrix translates immediately into a relation satisfied by our
family of integrals.  The ideal of such relations is known to be generated
by a family of bilinear equations, known as the Pl\"ucker relations.  For
our purposes, we restrict our attention to the simplest such identities.

For any matrix $M$, we let $\det_{S,T}(M)$ denote the determinant of the
submatrix of $M$ with coordinates $i\in S,j\in T$.  (This has a sign
ambiguity which we can eliminate by fixing an ordering on the coordinates.)
These determinants satisfy the following three identities:
\[
\det_{S\cup \{a,b\},T\cup\{c,d\}}(M)\det_{S,T}(M)
-
\det_{S\cup\{a\},T\cup\{c\}}(M)\det_{S\cup\{b\},T\cup\{d\}}(M)
+
\det_{S\cup\{a\},T\cup\{d\}}(M)\det_{S\cup\{b\},T\cup\{c\}}(M)
=
0,
\]
where $|S|=|T|$ and the coordinates are ordered so that $S<a<b;T<c<d$,
\[
\det_{S\cup\{a\},T\cup\{c,d\}}(M)
\det_{S,T\cup\{b\}}(M)
-
\det_{S\cup\{a\},T\cup\{b,d\}}(M)
\det_{S,T\cup\{c\}}(M)
+
\det_{S\cup\{a\},T\cup\{b,c\}}(M)
\det_{S,T\cup\{d\}}(M)
=
0,
\]
where $|S|=|T|+1$ and the coordinates are ordered so that $S<a;T<b<c<d$,
and
\[
\det_{S,T\cup\{a,b\}}(M)
\det_{S,T\cup\{c,d\}}(M)
-
\det_{S,T\cup\{a,c\}}(M)
\det_{S,T\cup\{b,d\}}(M)
+
\det_{S,T\cup\{a,d\}}(M)
\det_{S,T\cup\{b,c\}}(M)
=
0,
\]
where $|S|=|T|+2$ and the coordinates are ordered so that $T<a<b<c<d$.

Applying these identities to the matrix with entries
\[
\tau^{(1)}(\mu[a_i][b_j]'),
\]
and rescaling $\mu$, we obtain the following identities.

\begin{thm}\label{thm:detFay}
Assuming the integrals in question are all defined, we have the following
identities.
\begin{align}
 \tau^{(n+1)}(\mu[a,b][c,d]')\tau^{(n-1)}(\mu)
-\tau^{(n)}(\mu[a][c]')\tau^{(n)}(\mu[b][d]')+
\tau^{(n)}(\mu[a][d]')\tau^{(n)}(\mu[b][c]')&=0\\
  \tau^{(n-1)}(\mu[b])\tau^{(n)}(\mu[c,d][a]')
-\tau^{(n-1)}(\mu[c])\tau^{(n)}(\mu[b,d][a]')
 +\tau^{(n-1)}(\mu[d])\tau^{(n)}(\mu[c,d][a]')&=0\\
 \tau^{(n)}(\mu[c,d])\tau^{(n)}(\mu[a,b])
-\tau^{(n)}(\mu[b,d])\tau^{(n)}(\mu[a,c])
+\tau^{(n)}(\mu[b,c])\tau^{(n)}(\mu[a,d])&=0.
\end{align}
If we set $\tau^{(n)}=0$ for $n<0$, these identities remain valid for all
integers $n$.
\end{thm}

\begin{proof}
The Pl\"ucker identity argument immediately gives the first identity for
$n\ge 1$ and the other identities for $n\ge 2$.  Similarly, the first two
identities are trivial for $n\le 0$, and the third identity is trivial for
$n\le -1$.  So it remains to show the second identity for $n=1$ and the
third identity for $n=0$, $n=1$.  The second identity for $n=1$ is a linear
relation between univariate integrals that follows immediately from the
relation
\[
\psi(c,d)\psi(b,x)-\psi(b,d)\psi(d,x)+\psi(b,c)\psi(c,x).
\]
The third identity for $n=0$ is just the case $x=a$ of this identity.
Finally, for $n=1$, the third identity is the pfaffian of a $4\times 4$
matrix which has rank 2 by the second identity.
\end{proof}

\begin{rem}
When $\psi(x,y)=\psi'(x,y)=x-y$, these are instances of the generalized Fay
identities of \cite{AdlerM/ShiotaT/vanMoerbekeP:1998}.  By the remark after
Proposition \ref{prop:ell_CauchyII}, this can be used to obtain the general
$\psi=\psi'$ case via a change of variables.  Similarly, the case
$\psi\ne\psi'$ can be obtained via a change of variables and a delta
function limit from the identities of \cite{AdlerM/vanMoerbekeP:1999}.  
The above more elementary proof based on the Cauchy determinant appears to
be new, however.
\end{rem}

Similarly, if $\epsilon$ is an arbitrary antisymmetric function on $X$,
define, for $n$ even,
\[
\tau^{(n)}_{1/2}(\mu;\epsilon;\psi)
=
\frac{1}{n!}
\int
\pf_{1\le i,j\le n}(\epsilon(x_i,x_j))
\prod_{1\le i<j\le n} \psi(x_i,x_j)
\prod_{1\le i\le n} \mu(dx_i)
\]
For $n$ odd, we also need a univariate function $\phi$, and then define
\[
\tau^{(n)}_{1/2}(\mu;\phi,\epsilon;\psi)
=
\frac{1}{n!}
\int
\pf_{1\le i,j\le n}(\phi(x_i);\epsilon(x_i,x_j))
\prod_{1\le i<j\le n} \psi(x_i,x_j)
\prod_{1\le i\le n} \mu(dx_i)
\]
Here, for $n$ odd,
\[
\pf_{1\le i,j\le n}(\phi(x_i);\epsilon(x_i,x_j))
\]
represents the pfaffian of the $n+1\times n+1$ antisymmetric matrix
obtained from the $n\times n$ matrix $\epsilon(x_i,x_j)$ by adjoining a row
$\phi(x_i)$ and a column $-\phi(x_i)$.  Again, we will fix
$\phi,\epsilon,\psi$ and suppress them from the notation; whether $\phi$
appears is determined from the parity of $n$.  Similarly, we extend the
notation to cover negative integers by setting $\tau^{(n)}_{1/2}=0$ for
$n<0$.

It follows by an identity of de Bruijn \cite{deBruijnNG:1955} (also see
\cite{TracyCA/WidomH:1998} for a discussion in the context of Selberg
integrals) that these pfaffian $\tau$ functions can be written as
pfaffians.

\begin{prop}
Assuming the integrals are all defined, we have the following expressions.
For $n$ even,
\begin{align}
\tau^{(n)}_{1/2}(\mu[a_1,\dots,a_n])
&=
\pf_{1\le i,j\le n}\left(\tau^{(2)}_{1/2}(\mu[a_i,a_j])\right)\\
\noalign{\noindent and for $n$ odd,}
\tau^{(n)}_{1/2}(\mu[a_1,\dots,a_n])
&=
\pf_{1\le i,j\le n}\left(
\tau^{(1)}_{1/2}(\mu[a_i])
;
\tau^{(2)}_{1/2}(\mu[a_i,a_j])\right)
\end{align}
\end{prop}

Similarly to the determinantal case, there are a number of bilinear
identities satisfied by the pfaffian minors of an antisymmetric matrix.
For our purposes, we will restrict our attention to the following pair of
four-term identities.
\begin{align}
 \pf_S(A)\pf_{S\cup\{a,b,c,d\}}(A)
-\pf_{S\cup\{a,b\}}(A)\pf_{S\cup\{c,d\}}(A)&\notag\\
{}+\pf_{S\cup\{a,c\}}(A)\pf_{S\cup\{b,d\}}(A)
-\pf_{S\cup\{a,d\}}(A)\pf_{S\cup\{b,c\}}(A)
&=
0,
\end{align}
where $|S|$ is even, $S<a<b<c<d$, and
\begin{align}
 \pf_{S\cup\{a\}}(A)\pf_{S\cup\{b,c,d\}}(A)
-\pf_{S\cup\{b\}}(A)\pf_{S\cup\{a,c,d\}}(A)&\notag\\
{}+\pf_{S\cup\{c\}}(A)\pf_{S\cup\{a,b,d\}}(A)
-\pf_{S\cup\{d\}}(A)\pf_{S\cup\{a,b,c\}}(A)
&=
0,
\end{align}
where $|S|$ is odd, $S<a<b<c<d$.  See \cite{KnuthDE:1996} for these, and other
such identities.

These give rise to identities between our pfaffian $\tau$ functions; since in
the odd-dimensional case one row and column is special, we obtain a total
of six such identities.  However, the resulting identities turn out to
behave the same for $n$ odd and $n$ even, giving us a total of three
identities.

\begin{thm}
Assuming the integrals are all defined, we have the following identities
for all integers $n$.
\begin{align}
\tau^{(n+4)}_{1/2}(\mu[a_1,a_2,a_3,a_4])
\tau^{(n)}_{1/2}(\mu)
-
\tau^{(n+2)}_{1/2}(\mu[a_1,a_2])
\tau^{(n+2)}_{1/2}(\mu[a_3,a_4])&\notag\\
{}+
\tau^{(n+2)}_{1/2}(\mu[a_1,a_3])
\tau^{(n+2)}_{1/2}(\mu[a_2,a_4])
-
\tau^{(n+2)}_{1/2}(\mu[a_1,a_4])
\tau^{(n+2)}_{1/2}(\mu[a_2,a_3])
&{}=0.\\
\tau^{(n+3)}_{1/2}(\mu[a_2,a_3,a_4])
\tau^{(n+1)}_{1/2}(\mu[a_1])
-
\tau^{(n+3)}_{1/2}(\mu[a_1,a_3,a_4])
\tau^{(n+1)}_{1/2}(\mu[a_2])&\notag\\
{}+
\tau^{(n+3)}_{1/2}(\mu[a_1,a_2,a_4])
\tau^{(n+1)}_{1/2}(\mu[a_3])
-
\tau^{(n+3)}_{1/2}(\mu[a_1,a_2,a_3])
\tau^{(n+1)}_{1/2}(\mu[a_4])
&{}=0.\\
\tau^{(n+3)}_{1/2}(\mu[a_1,a_2,a_3])
\tau^{(n)}_{1/2}(\mu)
-
\tau^{(n+2)}_{1/2}(\mu[a_2,a_3])
\tau^{(n+1)}_{1/2}(\mu[a_1])
&\notag\\
{}+
\tau^{(n+2)}_{1/2}(\mu[a_1,a_3])
\tau^{(n+1)}_{1/2}(\mu[a_2])
-
\tau^{(n+2)}_{1/2}(\mu[a_1,a_2])
\tau^{(n+1)}_{1/2}(\mu[a_3])
&{}=0.
\end{align}
\end{thm}

\begin{proof}
For $n\ge 0$, these identities are just relations between minors of the
antisymmetric matrix
\[
(\tau^{(1)}_{1/2}(\mu[a_i]);\tau^{(2)}_{1/2}(\mu[a_i,a_j]))
\]
For $n\le -1$, the first and third identities follow from Theorem
\ref{thm:detFay} (as, when nontrivial, they relate univariate integrals and
scalars).  The only remaining nontrivial case is the instance $n=-1$ of the
second identity.  But this is a linear relation between bivariate integrals
coming from a linear relation of the integrands.
\end{proof}

\begin{rem}
Again, these could be obtained via a change of variables from the pfaffian
Fay identities of \cite{AdlerM/ShiotaT/vanMoerbekeP:2002}, but our
elementary proof is new.
\end{rem}

\section{Painlev\'e recurrences}

If we apply the generalized Fay identities to an integral of the form
$\tilde\II_{q;p,q}$, we find that for suitable choices of $a_i$, the
integrals that appear are of the same form, with shifted parameters.  As a
result, Theorem \ref{thm:detFay} gives rise to three special cases of the
following recurrence.

\begin{thm}
Let $v_0$, $v_1$, $v_2\in \frac{1}{2}\Lambda_{E_8}$ be unit vectors in a
common coset of $\Lambda_{E_8}$, and let $\phi\in \Hom(\Lambda_{E_8},\C^*)$
be such that $pq/(\tau_{v_0}\phi)(\omega)\in q^\Z$.  Then
\[
\sum_{0\le r\le 2}
\frac{
\tilde\II_{q;p,q}(\tau_{v_r}\phi)
\tilde\II_{q;p,q}(\tau_{-v_r}\phi)}
{
\prod_{s\ne r} \psi_p(\phi(v_r),\phi(v_s))
}
=
0.
\]
\end{thm}

\begin{proof}
First, we observe as remarked that Theorem \ref{thm:detFay} gives
essentially three special cases (up to signed permutations within the
triples and the natural action of $S_8$ on the coordinates), namely:
\begin{align}
v_0&=(\frac{1}{2},\frac{1}{2},\frac{1}{2},\frac{1}{2},0,0,0,0)&v_1&=(\frac{1}{2},\frac{1}{2},-\frac{1}{2},-\frac{1}{2},0,0,0,0)&v_2&=(\frac{1}{2},-\frac{1}{2},-\frac{1}{2},\frac{1}{2},0,0,0,0)\\
v_0&=(\frac{1}{2},-\frac{1}{2},\frac{1}{2},\frac{1}{2},0,0,0,0)&v_1&=(\frac{1}{2},\frac{1}{2},-\frac{1}{2},\frac{1}{2},0,0,0,0)&v_2&=(\frac{1}{2},\frac{1}{2},\frac{1}{2},-\frac{1}{2},0,0,0,0)\\
v_0&=(\frac{1}{2},\frac{1}{2},-\frac{1}{2},-\frac{1}{2},0,0,0,0)&v_1&=(\frac{1}{2},-\frac{1}{2},\frac{1}{2},-\frac{1}{2},0,0,0,0)&v_2&=(\frac{1}{2},-\frac{1}{2},-\frac{1}{2},\frac{1}{2},0,0,0,0)
\end{align}

Now, we know that $\tilde\II_{q;p,q}$ is invariant under the action of
$W(E_7)$, so it suffices to prove the theorem for one (unsigned, unordered)
triple from each $W(E_7)$-orbit.  There are four such orbits, so we still
have one orbit remaining to consider, one representative of which is:
\[
v_0=(\frac{1}{2},\frac{1}{2},\frac{1}{2},\frac{1}{2},0,0,0,0)\quad
v_1=(0,0,0,0,\frac{1}{2},\frac{1}{2},\frac{1}{2},\frac{1}{2})\quad
v_2=(0,0,0,0,\frac{1}{2},\frac{1}{2},-\frac{1}{2},-\frac{1}{2}).
\]
Now, consider the following two representatives of the orbit of our first
special case:
\begin{align}
v_0&=(\frac{1}{2},\frac{1}{2},\frac{1}{2},\frac{1}{2},0,0,0,0)&
v_1&=(0,0,0,0,\frac{1}{2},-\frac{1}{2},\frac{1}{2},-\frac{1}{2})&
v_2&=(0,0,0,0,\frac{1}{2},\frac{1}{2},-\frac{1}{2},-\frac{1}{2})\\
v_0&=(0,0,0,0,\frac{1}{2},\frac{1}{2},\frac{1}{2},\frac{1}{2})&
v_1&=(0,0,0,0,\frac{1}{2},-\frac{1}{2},\frac{1}{2},-\frac{1}{2})&
v_2&=(0,0,0,0,\frac{1}{2},\frac{1}{2},-\frac{1}{2},-\frac{1}{2}).
\end{align}
If we take a linear combination of these two bilinear identities in such
a way as to eliminate the term corresponding to
$(0,0,0,0,1/2,-1/2,1/2,-1/2)$, we find that the result is precisely the
desired bilinear identity corresponding to the above representative of
the missing orbit.
\end{proof}

As the above scheme of bilinear recurrences has appeared elsewhere in the
literature \cite[Theorem 5.2]{KajiwaraK/MasudaT/NoumiM/OhtaY/YamadaY:2004},
we immediately obtain the following corollary.  Compare also the results of
\cite{ForresterPJ/WitteNS:2002} for the Selberg limit.

\begin{cor}
The function $\tilde\II_{q;p,q}$ is a tau function for the elliptic
Painlev\'e equation.
\end{cor}

\begin{rem}
We should mention in this context that Sakai's version of elliptic
Painlev\'e \cite{SakaiH:2001} is geometrically described in terms of the
blow-up of $\P^2$ at 9 points.  This has a natural $S_9$ symmetry, but the
$S_8$ symmetry it gives is {\em not} conjugate to the natural $S_8$
symmetry on our integral.  This suggests that, at least from the integral
perspective, the more natural geometric context for elliptic Painlev\'e is
the blow-up of $\P^1\times \P^1$ at 8 points.
\end{rem}

We next turn to the case $t=q^2$.  Here, the cross term in the integral is
\[
\prod_{1\le i<j\le n}
\theta_p(z_i^{\pm 1}z_j^{\pm 1})
\theta_p(q z_i^{\pm 1}z_j^{\pm 1}).
\]
As it stands, this does not appear to be amenable to either generalized
Fay identity.  However, we can write this as
\[
q^{n(n-1)}
\prod_{1\le i<j\le n}
\psi_p(q^{\pm 1/2} z_i,q^{\pm 1/2} z_j),
\]
which can in turn be written as
\[
q^{5n(n-1)/4}
\prod_{1\le i\le n} \frac{1}{z_i^{-1}\theta_p(z_i^2)}
\prod_{1\le i<j\le 2n} \psi_p(w_i,w_j),
\]
where $w_{2i-1} = q^{1/2}z_i$, $w_{2i}=q^{-1/2}z_i$.  But, aside from a
factor of $(2n)!/2^n n!$, this corresponds to a limiting case of
$\tau^{(2n)}_{1/2}$, taking
\[
\epsilon(x_i,x_j) = \delta(x_i-qx_j)-\delta(x_j-qx_i).
\]
As a result, we again obtain bilinear identities, this time with four
terms each.

\begin{thm}
Let $v_0$, $v_1$, $v_2$, $v_3\in \frac{1}{2}\Lambda_{E_8}$ be unit vectors
in a common coset of $\Lambda_{E_8}$, and let $\phi\in
\Hom(\Lambda_{E_8},\C^*)$ such that $pq/(\tau_{2v_0}\phi)(\omega)\in
q^{2\Z}$.  Then
\[
\sum_{0\le r\le 3}
\frac{
\tilde\II_{q^2;p,q}(\tau_{2v_r}\phi)
\tilde\II_{q^2;p,q}(\tau_{-2v_r}\phi)}
{
\prod_{s\ne r} \psi_p(\phi(v_r),\phi(v_s))
}
=
0.
\]
\end{thm}

\begin{proof}
Here, the generalized Fay identities give us two of the four
$W(E_7)$-orbits we require:
\begin{align}
v_0&=(\frac{1}{2},\frac{1}{2},\frac{1}{2},\frac{1}{2},0,0,0,0),&
v_1&=(\frac{1}{2},\frac{1}{2},-\frac{1}{2},-\frac{1}{2},0,0,0,0),\notag\\
v_2&=(\frac{1}{2},-\frac{1}{2},\frac{1}{2},-\frac{1}{2},0,0,0,0),&
v_3&=(\frac{1}{2},-\frac{1}{2},-\frac{1}{2},\frac{1}{2},0,0,0,0),\\
\noalign{\noindent and}
v_0&=(-\frac{1}{2},\frac{1}{2},\frac{1}{2},\frac{1}{2},0,0,0,0),&
v_1&=(\frac{1}{2},-\frac{1}{2},\frac{1}{2},\frac{1}{2},0,0,0,0),\notag\\
v_2&=(\frac{1}{2},\frac{1}{2},-\frac{1}{2},\frac{1}{2},0,0,0,0),&
v_3&=(\frac{1}{2},\frac{1}{2},\frac{1}{2},-\frac{1}{2},0,0,0,0).
\end{align}
But again, we can obtain the bilinear identities corresponding to the
missing orbits as linear combinations of $W(E_7)$-images of the identities
coming from the generalized Fay identities.  The key point here is that
if two of our identities have three monomials in common, and we take a
linear combination to eliminate one of those monomials, the result is
another of our identities.
\end{proof}

For $t=q^{1/2}$, we apply Theorem \ref{thm:ell_CauchypfaffII} (as an
identity of $q$-theta functions!) to write the integrand as an instance of
$\tau^{(n)}_{1/2}$.  The first two generalized Fay identities turn out to
give us a recurrence strikingly like the one we have just seen; again, by
taking linear combinations of $W(E_7)$-images, we can obtain the entire
$W(E_8)$-orbit.

\begin{thm}
Let $v_0$, $v_1$, $v_2$, $v_3\in \frac{1}{2}\Lambda_{E_8}$ be unit vectors in a
common coset of $\Lambda_{E_8}$, and let
$\phi\in \Hom(\Lambda_{E_8},\C^*)$ such that
$pq/(\tau_{v_0}\phi)(\omega)\in q^{\Z/2}$.
Then
\[
\sum_{0\le r\le 3}
\frac{
\tilde\II_{q^{1/2};p,q}(\tau_{v_r}\phi)
\tilde\II_{q^{1/2};p,q}(\tau_{-v_r}\phi)}
{
\prod_{s\ne r} \psi_p(\phi(v_r),\phi(v_s))
}
=
0.\label{eq:EP_beta_one}
\]
\end{thm}

The third Fay identity corresponds to a different $W(E_8)$-orbit.
Unfortunately, we can no longer combine instances of the $W(E_7)$-orbit,
and have in fact been unable to prove the presumable general form of the
recurrence.  We do, however, have the following.

\begin{thm}
Suppose the unordered, unsigned quadruple $(v_0,v_1,v_2,v_3)$ is in the
$W(E_7)$-orbit of
\[
(
(\frac{1}{2},\frac{1}{2},\frac{1}{2},0,0,0,0,0),
(\frac{1}{2},-\frac{1}{2},-\frac{1}{2},0,0,0,0,0),
(-\frac{1}{2},-\frac{1}{2},\frac{1}{2},0,0,0,0,0),
(-\frac{1}{2},\frac{1}{2},-\frac{1}{2},0,0,0,0,0)
).
\]
Then equation \eqref{eq:EP_beta_one} still holds.
\end{thm}

Unlike the determinantal case, the recurrences corresponding to the
pfaffian cases ($t\in\{q^2,q^{1/2}\}$) appear to be new, even in the
Painlev\'e setting.  It would be very interesting to know a geometric
interpretation for these recurrences, and more generally to understand how
they relate to the usual elliptic Painlev\'e equation.

\end{document}